\documentclass[a4paper, 12pt]{article}
\usepackage{amsfonts}
\usepackage{fancyhdr}
\usepackage{graphicx}
\usepackage{epsfig}
\usepackage{amssymb}
 \usepackage{mathrsfs,amsfonts,amsmath}
 \usepackage{color}

\makeatletter

\@addtoreset{equation}{section} \makeatother
\newtheorem{Theorem}{Theorem}[section]
\newtheorem{Remark}[Theorem]{Remark}
\newtheorem{Lemma}[Theorem]{Lemma}

\begin{document}
\title{\textbf{Exponential stability of modified truncated EM method for stochastic differential equations
\footnote{Supported by Natural Science Foundation of China (NSFC
11601025).}}}
\author{ Guangqiang Lan\footnote{Corresponding author: Email:
langq@mail.buct.edu.cn.}\quad and\quad
 Fang Xia
\\ \small School of Science, Beijing University of Chemical Technology, Beijing 100029, China}

\date{}

\maketitle

\begin{abstract}
Exponential stability of modified truncated Euler-Maruyama method for stochastic differential equations are investigated in this paper. Sufficient conditions for the $p$-th moment and almost sure exponential stability of the given numerical method are presented. An example is provided to support our conclusions.
\end{abstract}



\section{Introduction and main results}

\noindent

Stability theorems of stochastic differential systems, for example, moment stability (M-stability) and almost sure stability (or the
trajectory stability (T-stability)), have attracted more and more attention in recent years. Since the exact solution is usually difficult to obtain, properties of different types of numerical simulations are more and more hot topics. There are many results on asymptotic stability theorems for different stochastic numerical approximations. For example,
\cite{BB,HMY,LFM,MM} considered the exponential (or polynomial) stability of the corresponding EM method and \cite{HMS,LFM,MS1,WMS} considered backward EM method. There also many results for more general stochastic theta method (it is also called $\theta$-EM method), such as \cite{CW,HI,LY,MS,RS,WC,ZW} etc. However, the classical EM method usually needs the linear growth condition for both coefficients and the semi-implicit $\theta$-EM method ($0<\theta\le1$) needs the one-sided Lipschitz condition for the drift term $f$ to ensure the stability of the corresponding numerical method. So we will investigate the stability of a new numerical method (we call it modified truncated Euler-Maruyama method) which does need such conditions.

Recently, Professor Mao introduced in \cite{mao} a new numerical simulation called truncated EM method, and then he obtained sufficient conditions for the strong convergence rate of it in \cite{mao1}. However, to the best of our knowledge, stability of truncated EM method is still unknown. In this work, we
will study the exponential stability of the simulations of the given stochastic differential equations.

Let $(\Omega,\mathscr{F},(\mathscr{F}_t)_{t\geq 0},P)$ be a complete
filtered probability space satisfying usual conditions. Consider the following stochastic differential
equations:

\begin{equation}\label{sde}dx(t)
=f(x(t))dt+g(x(t))dB_t,\ X_0=x_0\in\mathbb{R}^d
\end{equation}
where  $(B_t)_{t\geq0}$ is an standard scalar
$\mathscr{F}_t$-Brownian motion,
$f:x\in\mathbb{R}^d\mapsto
f(x)\in\mathbb{R}^d$ and
$g:x\in\mathbb{R}^d\mapsto g(x)\in\mathbb{R}^d$ are measurable functions.

Suppose
$$f(0)\equiv0,\quad g(0)\equiv0,$$
which implies that $X\equiv0$ is the trivial solution of equation
(\ref{sde}).

Assume that the coefficients satisfy local Lipschitz condition, that is, for each
$R$ there is $L_R>0$ (depending on $R$) such that
\begin{equation}\label{local}|f(x)-f(\bar{x})|
\vee|g(x)-g(\bar{x})|\le L_R|x-\bar{x}|\end{equation} for all
$|x|\vee|\bar{x}|\le R>0$.

It is obvious that $L_R$ is increasing function with respect to $R$. It is also well known that there is a unique strong solution (might explode at finite time) to equation (\ref{sde}) under local Lipschitz condition (\ref{local}) (Indeed, local Lipschitz condition could be relaxed to non-Lipschitz condition, see e.g. \cite{LW}).

Choose $\Delta^*>0$ small enough and a strictly positive decreasing function $h:(0,\Delta^*]\to(0,\infty)$ such that
\begin{equation}\label{tiaoj}\lim_{\Delta\to0}h(\Delta)=\infty\ \textrm{and}\ \lim_{\Delta\to0}L_{h(\Delta)}^4\Delta=0.\end{equation}

\begin{Remark}
Notice that such function $h$ always exists for given Lipschitz coefficient $L_R$. For example, let $l(R)=\frac{1}{RL_R^4}$ and $h$ is the inverse function of $l$. Then $h$ is decreasing and $\lim_{\Delta\to0}h(\Delta)=\infty$ since $l$ is decreasing and $\lim_{R\to\infty}l(R)=0$. If we set $R=h(\Delta),$ then $L_{h(\Delta)}^4\Delta=L_R^4l(R)=\frac{1}{R}=\frac{1}{h(\Delta)}\to0$ as $\Delta\to0$.
\end{Remark}

Motivated by Mao \cite{mao}, for any $\Delta>0,$ we define the modified truncated function of $f$ as the following:
$$f_\Delta(x)=\left\{\begin{array}{ll} f(x),|x|\le h(\Delta),\\
\frac{|x|}{h(\Delta)} \cdot f(h(\Delta)\cdot\frac{x}{|x|}),|x|>h(\Delta).
\end{array}
\right.$$
$g_\Delta$ is defined in the same way as $f_\Delta$.

Notice that the function $f_\Delta$ defined above is different from Mao \cite{mao,mao1} (where the truncated functions are bound for any fixed $\Delta$).

Then we define the modified truncated EM method (MTEM) numerical solutions $X_k^\Delta\approx x(k\Delta)$ by setting $X_0^\Delta=x_0$ and
\begin{equation}\label{num}X_{k+1}^\Delta=X_k^\Delta+f_\Delta(X_k^\Delta)\Delta+g_\Delta(X_k^\Delta)\Delta B_k\end{equation}
for $k=0,1,2,\cdots,$ where $\Delta B_k=B((k+1)\Delta)-B(k\Delta)$ is the increment of the scalar Brownian motion.

The two versions of the continuous-time MTEM solutions are defined as the following:

\begin{equation}\label{num1}\bar{x}_\Delta(t)=\sum_{k=0}^\infty X_k^\Delta1_{[k\Delta,(k+1)\Delta)}(t),\quad t\ge0,\end{equation}
and
\begin{equation}\label{num2}x_\Delta(t)=x_0+\int_0^tf_\Delta(\bar{x}_\Delta(s))ds+\int_0^tg_\Delta(\bar{x}_\Delta(s))dB(s),\quad t\ge0.\end{equation}

It is easy to see that $x_\Delta(k\Delta)=\bar{x}_\Delta(k\Delta)=X_k^\Delta$ for all $k\ge0.$

In \cite{LX}, the authors investigated the strong convergence of the two version of continuous-time MTEM approximations (\ref{num1}) and (\ref{num2}). Now let us consider the exponential stability of the MTEM approximation (\ref{num}).

We have

\begin{Theorem}\label{stab}
Assume that $f$ and $g$ satisfy the local Lipschitz condition (\ref{local}),
and there exists sufficiently small $p\in(0,1)$ such that

\begin{equation}\label{tiaoj1}-\lambda:=\sup_{x\neq0}\left(\frac{\langle x,f(x)\rangle+\frac{1}{2}|g(x)|^2}{|x|^2}+\frac{p-2}{2}\frac{\langle x,g(x)\rangle^2}{|x|^4}\right)<0.\end{equation}

Then the solution of equation (\ref{sde}) satisfies
\begin{equation}\label{jingque}\limsup_{t\to\infty}\frac{\log\mathbb{E}(|x(t)|^p)}{t}\le-p\lambda.\end{equation}

If moreover, the local Lipschitz constant $L_R$ satisfies (\ref{tiaoj}), then for any sufficiently small $\varepsilon\in(0,1)$, there exists $\Delta^*>0$ such that for any $\Delta\in(0,\Delta^*)$ the discrete MTEM approximation (\ref{num}) satisfies
\begin{equation}\label{shuzhi}\limsup_{k\to\infty}\frac{\log\mathbb{E}(|X_{k}^\Delta|^p)}{k\Delta}\le-p(\lambda-\varepsilon).\end{equation}

Further, $X_{k}^\Delta$ is also almost surely exponentially stable of order
$\lambda-\varepsilon$, i.e.,
\begin{equation}\label{shuzhi1}\limsup_{k\to\infty}\frac{\log|X_{k}^\Delta|}{k\Delta}\le-(\lambda-\varepsilon),\quad a.s.\end{equation}

For the continuous-times MTEM approximation (\ref{num2}), we also have
\begin{equation}\label{shuzhi2}\limsup_{t\to\infty}\frac{\log\mathbb{E}(|x_\Delta(t)|^p)}{t}\le-p(\lambda-\varepsilon).\end{equation}

\end{Theorem}

\begin{Remark}
Under conditions (\ref{local}) and (\ref{tiaoj1}), equation (\ref{sde}) admits a unique strong global solution. Indeed, it is well known that (\ref{local}) indicates that the existence of a unique strong solution (which might explode at finite time) while (\ref{tiaoj1}) ensure the non explosion of the solution. Notice that since $0<p<1$,
$$\frac{\langle x,f(x)\rangle+\frac{1}{2}|g(x)|^2}{|x|^2}+\frac{p-2}{2}\frac{\langle x,g(x)\rangle^2}{|x|^4}\ge\frac{\langle x,f(x)\rangle+\frac{p-1}{2}|g(x)|^2}{|x|^2}.$$
Thus condition (\ref{tiaoj1}) implies the Khasminskii-type condition
$$\langle x,f(x)\rangle+\frac{p-1}{2}|g(x)|^2\le K(1+|x|^2).$$
\end{Remark}

The organization of the paper is as the following. In Section 2, global Lipshitz continuity of the modified truncated functions $f_\Delta$ and $g_\Delta$ is proved. In Section 3, the exponential stability of MTEM is proved. Finally, an example is presented to interpret the theory.

\section{Global Lipshitz continuity of $f_\Delta$ and $g_\Delta$}

For the modified truncated function $f_\Delta$ and $g_\Delta$, we have the following

\begin{Lemma}\label{l1}
Suppose the local Lipschitz condition (\ref{local}) holds. Then for any fixed $\Delta>0$,
\begin{equation}\label{global}|f_\Delta(x)-f_\Delta(\bar{x})|
\vee|g_\Delta(x)-g_\Delta(\bar{x})|\le 3L_{h(\Delta)}|x-\bar{x}|,\ \forall x,\bar{x}\in\mathbb{R}^d.\end{equation}
\end{Lemma}

\textbf{Proof}\ For any $x,\bar{x}\in\mathbb{R}^d$, there are three cases: $x,\bar{x}$ are both in the ball $B(h(\Delta))=\{x\in\mathbb{R}^d,|x|\le h(\Delta)\},$ $x,\bar{x}$ are both outside the ball $B(h(\Delta))$ and one is in the ball and the other is outside the ball.

If $x,\bar{x}\le h(\Delta)$, then (\ref{global}) holds naturally by (\ref{local});
Now assume $x,\bar{x}> h(\Delta)$. Since
$$|h(\Delta)\cdot\frac{x}{|x|}|=|h(\Delta)\cdot\frac{\bar{x}}{|\bar{x}|}|=h(\Delta),$$
then we have
$$\aligned |f_\Delta(x)-f_\Delta(\bar{x})|&=\left|\frac{|x|}{h(\Delta)} \cdot f\left(h(\Delta)\cdot\frac{x}{|x|}\right)-\frac{|\bar{x}|}{h(\Delta)} \cdot f\left(h(\Delta)\cdot\frac{\bar{x}}{|\bar{x}|}\right)\right|\\&
\le \frac{|x|}{h(\Delta)}\left|f\left(h(\Delta)\cdot\frac{x}{|x|}\right)-f\left(h(\Delta)\cdot\frac{\bar{x}}{|\bar{x}|}\right)\right|\\&
\quad+ \frac{\left||x|-|\bar{x}|\right|}{h(\Delta)}\left|f\left(h(\Delta)\cdot\frac{\bar{x}}{|\bar{x}|}\right)\right|\\&
\le \frac{|x|}{h(\Delta)}\cdot L_{h(\Delta)}\left|h(\Delta)\cdot\frac{x}{|x|}-h(\Delta)\cdot\frac{\bar{x}}{|\bar{x}|}\right|\\&
\quad+ \frac{\left||x|-|\bar{x}|\right|}{h(\Delta)}\left(L_{h(\Delta)}\left|h(\Delta)\cdot\frac{\bar{x}}{|\bar{x}|}\right|\right)\\&
\le L_{h(\Delta)}\left|x-\frac{|x|}{|\bar{x}|}\bar{x}\right|+L_{h(\Delta)}|x-\bar{x}|\\&
\le L_{h(\Delta)}\left(\left|x-\bar{x}\right|+\left|\bar{x}-\frac{|x|}{|\bar{x}|}\bar{x}\right|\right)
+L_{h(\Delta)}|x-\bar{x}|\\&
\le 3L_{h(\Delta)}|x-\bar{x}|.\endaligned$$

Finally, without loss of generality, suppose that $|x|\le h(\Delta)<|\bar{x}|.$ Then we have
$$\aligned |f_\Delta(x)-f_\Delta(\bar{x})|&=\left|f(x)-\frac{|\bar{x}|}{h(\Delta)} \cdot f\left(h(\Delta)\cdot\frac{\bar{x}}{|\bar{x}|}\right)\right|\\&
\le \left|f(x)-f\left(h(\Delta)\cdot\frac{\bar{x}}{|\bar{x}|}\right)\right|
+\left|f\left(h(\Delta)\cdot\frac{\bar{x}}{|\bar{x}|}\right)\right|\left|1-\frac{|\bar{x}|}{h(\Delta)}\right|\\&
\le L_{h(\Delta)}\left|x-h(\Delta)\cdot\frac{\bar{x}}{|\bar{x}|}\right|+L_{h(\Delta)}|h(\Delta)-|\bar{x}||.
\endaligned$$

Since $|x|\le h(\Delta)<|\bar{x}|$, then $|h(\Delta)-|\bar{x}||=|\bar{x}|-h(\Delta)\le|\bar{x}|-|x|\le|x-\bar{x}|$, and $$\aligned\left|x-h(\Delta)\cdot\frac{\bar{x}}{|\bar{x}|}\right|&\le |x-\bar{x}|+\left|\bar{x}-h(\Delta)\cdot\frac{\bar{x}}{|\bar{x}|}\right|\\&
=|x-\bar{x}|+|h(\Delta)-|\bar{x}||\le2|x-\bar{x}|.\endaligned$$

Therefore,
$$\aligned|f_\Delta(x)-f_\Delta(\bar{x})|&\le 3L_{h(\Delta)}|x-\bar{x}|.\endaligned$$

Similarly, $g_\Delta$ is also globally Lipschitz continuous with the same Lipschitz constant $3L_{h(\Delta)}$. We complete the proof. $\square$

\begin{Lemma}\label{l2}
Suppose (\ref{tiaoj1}) holds for some $p\in(0,1)$. Then we have
\begin{equation}\label{tiaoj2}\sup_{x\neq0}\left(\frac{\langle x,f_\Delta(x)\rangle+\frac{1}{2}|g_\Delta(x)|^2}{|x|^2}+\frac{p-2}{2}\frac{\langle x,g_\Delta(x)\rangle^2}{|x|^4}\right)\le-\lambda.\end{equation}
\end{Lemma}

\textbf{Proof}\ By definition of the modified truncated functions $f_\Delta$ and $g_\Delta,$ we have
$$\aligned &\quad\sup_{x\neq0}\left(\frac{\langle x, f_\Delta(x)\rangle+\frac{1}{2}|g_\Delta(x)|^2}{|x|^2}+\frac{p-2}{2}\frac{\langle x, g_\Delta(x)\rangle^2}{|x|^4}\right)\\&
=\max\left\{\sup_{0<|x|\le h(\Delta)}\left(\frac{\langle x, f(x)\rangle+\frac{1}{2}|g(x)|^2}{|x|^2}+\frac{p-2}{2}\frac{\langle x, g(x)\rangle^2}{|x|^4}\right),\right.\\&
\qquad\qquad\sup_{|x|>h(\Delta)}\left(\frac{\langle x, \frac{|x|}{h(\Delta)} \cdot f(h(\Delta)\cdot\frac{x}{|x|})\rangle+\frac{1}{2}|\frac{|x|}{h(\Delta)} \cdot g(h(\Delta)\cdot\frac{x}{|x|})|^2}{|x|^2}\right.\\&\left.\left.\qquad\qquad\qquad\qquad+\frac{p-2}{2}\frac{\langle x, \frac{|x|}{h(\Delta)} \cdot g(h(\Delta)\cdot\frac{x}{|x|})\rangle^2}{|x|^4}\right)\right\}.\endaligned$$
If we set $y=h(\Delta)\cdot\frac{x}{|x|}$, then $|y|=h(\Delta)$, and therefore
$$\aligned &\quad\sup_{x\neq0}\left(\frac{\langle x, f_\Delta(x)\rangle+\frac{1}{2}|g_\Delta(x)|^2}{|x|^2}+\frac{p-2}{2}\frac{\langle x, g_\Delta(x)\rangle^2}{|x|^4}\right)\\&
=\max\left\{\sup_{0<|x|\le h(\Delta)}\left(\frac{\langle x, f(x)\rangle+\frac{1}{2}|g(x)|^2}{|x|^2}+\frac{p-2}{2}\frac{\langle x, g(x)\rangle^2}{|x|^4}\right),\right.\\&
\left.\qquad\qquad\sup_{|x|>h(\Delta)}\left(\frac{\langle y, f(y)\rangle+\frac{1}{2}|g(y)|^2}{|y|^2}+\frac{p-2}{2}\frac{\langle y, g(y)\rangle^2}{|y|^4}\right)\right\}\\&
=\sup_{0<|x|\le h(\Delta)}\left(\frac{\langle x, f_\Delta(x)\rangle+\frac{1}{2}|g_\Delta(x)|^2}{|x|^2}+\frac{p-2}{2}\frac{\langle x, g_\Delta(x)\rangle^2}{|x|^4}\right)\le -\lambda.\endaligned$$

We complete the proof. $\square$

\section{Proof of Theorem \ref{stab}}
Firstly, let us prove the moment exponential stability of the exact solution $x(t).$ By It\^o's formula, we have
$$\aligned d(e^{p\lambda t}|x(t)|^p)&=e^{p\lambda t}|x(t)|^p\left[p\lambda+\frac{p}{2}\left(\frac{2\langle x(t), f(x(t))\rangle+|g(x(t))|^2}{|x(t)|^2}\right.\right.\\&
\qquad\qquad\qquad\left.\left.+(p-2)\frac{\langle x(t), g(x(t))\rangle^2}{|x(t)|^4}\right)\right]dt+dM(t),\endaligned$$
where
$$dM(t)=pe^{p\lambda t}|x(t)|^p\frac{\langle x(t), g(x(t))\rangle^2}{|x(t)|^2}dB(t)$$
is the martingale term.

Then by (\ref{tiaoj1}), we have
$$d(e^{p\lambda t}|x(t)|^p)\le dM(t).$$

Therefore,
$$e^{p\lambda t}\mathbb{E}(|x(t)|^p)\le |x_0|^p,$$
that is, the exact solution $x(t)$ is $p$-th moment exponentially stable.

Now let us prove (\ref{shuzhi}).

By definition of (\ref{num}), for any $k\ge0,$ we have

$$\aligned|X_{k+1}^\Delta|^2&=|X_{k}^\Delta|^2+2\langle X_{k}^\Delta,f_\Delta(X_k^\Delta)\Delta+g_\Delta(X_k^\Delta)\Delta B_k\rangle+|f_\Delta(X_k^\Delta)\Delta+g_\Delta(X_k^\Delta)\Delta B_k|^2\\&=|X_{k}^\Delta|^2(1+\xi_k),\endaligned$$
where
$$\xi_k=\frac{1}{|X_{k}^\Delta|^2}(2\langle X_{k}^\Delta,f_\Delta(X_k^\Delta)\Delta+g_\Delta(X_k^\Delta)\Delta B_k\rangle+|f_\Delta(X_k^\Delta)\Delta+g_\Delta(X_k^\Delta)\Delta B_k|^2)$$
if $X_{k}^\Delta\neq0$, otherwise, it is set to $-1.$ So we only need to suppose that $X_{k}^\Delta\neq0$.

Then we have
$$\aligned|X_{k+1}^\Delta|^p&=|X_{k}^\Delta|^p(1+\xi_k)^\frac{p}{2}.\endaligned$$

Notice that for any $x\ge-1$ and $p\in(0,1)$, by Taylor's expansion,
$$\aligned(1+x)^\frac{p}{2}&=1+\frac{p}{2}x+\frac{p(p-2)}{8}x^2+\frac{p(p-2)(p-4)}{2^3\times3!}x^3
+\frac{p(p-2)(p-4)(p-6)}{2^4\times4!}\theta^4
\\&\le1+\frac{p}{2}x+\frac{p(p-2)}{8}x^2+\frac{p(p-2)(p-4)}{2^3\times3!}x^3,\endaligned$$
where $\theta$ in the first equation lies between 0 and $x$.

Then we have
$$\mathbb{E}(|X_{k+1}^\Delta|^p|\mathscr{F}_{k\Delta})\le|X_{k}^\Delta|^p\mathbb{E}
(1+\frac{p}{2}\xi_k+\frac{p(p-2)}{8}\xi_k^2+\frac{p(p-2)(p-4)}{2^3\times3!}\xi_k^3|\mathscr{F}_{k\Delta}).$$

Now
$$\aligned\mathbb{E}(\xi_k|\mathscr{F}_{k\Delta})&=\frac{1}{|X_{k}^\Delta|^2}\mathbb{E}(2\langle X_{k}^\Delta,f_\Delta(X_k^\Delta)\Delta+g_\Delta(X_k^\Delta)\Delta B_k\rangle\\&\qquad\qquad+|f_\Delta(X_k^\Delta)\Delta+g_\Delta(X_k^\Delta)\Delta B_k|^2|\mathscr{F}_{k\Delta})\\&=\frac{1}{|X_{k}^\Delta|^2}[(2\langle X_{k}^\Delta, f_\Delta(X_k^\Delta)\rangle+|g_\Delta(X_k^\Delta)|^2)\Delta+|f_\Delta(X_k^\Delta)|^2\Delta^2].\endaligned$$

We have use the fact that $\mathbb{E}(\Delta B_k|\mathscr{F}_{k\Delta})=0$ and $\mathbb{E}((\Delta B_k)^2|\mathscr{F}_{k\Delta})=\Delta$ in the above equation.

Since $f_\Delta$ is globally Lipschitz continuous, then
$$\aligned\mathbb{E}(\xi_k|\mathscr{F}_{k\Delta})&\le\frac{2\langle X_{k}^\Delta, f_\Delta(X_k^\Delta)\rangle+|g_\Delta(X_k^\Delta)|^2}{|X_{k}^\Delta|^2}\Delta+L^2_{h(\Delta)}\Delta^2.\endaligned$$

Thus, by (\ref{tiaoj}), we have
$$\frac{L^2_{h(\Delta)}\Delta^2}{\Delta}=L^2_{h(\Delta)}\Delta\to0.$$

That is $L^2_{h(\Delta)}\Delta^2=o(\Delta).$ Here and from now on, $o(\Delta)$ represents the high order infinitesimal of $\Delta$ as $\Delta\to0.$ Therefore,
$$\aligned\mathbb{E}(\xi_k|\mathscr{F}_{k\Delta})&\le\frac{2\langle X_{k}^\Delta, f_\Delta(X_k^\Delta)\rangle+|g_\Delta(X_k^\Delta)|^2}{|X_{k}^\Delta|^2}\Delta+o(\Delta).\endaligned$$

Similarly, by (\ref{tiaoj}) and Lemma \ref{l1}, we have
$$\aligned\mathbb{E}(\xi^2_k|\mathscr{F}_{k\Delta})&=\frac{1}{|X_{k}^\Delta|^4}
\mathbb{E}((2\langle X_{k}^\Delta, g_\Delta(X_k^\Delta)\rangle\Delta B_k+B)^2|\mathscr{F}_{k\Delta})\\&=
\frac{1}{|X_{k}^\Delta|^4}\mathbb{E}(4\langle X_{k}^\Delta, g_\Delta(X_k^\Delta)\rangle^2(\Delta B_k)^2+B^2+4B\langle X_{k}^\Delta, g_\Delta(X_k^\Delta)\rangle\Delta B_k|\mathscr{F}_{k\Delta})\\&\ge \frac{4 \langle X_{k}^\Delta, g_\Delta(X_k^\Delta)\rangle^2}{|X_{k}^\Delta|^4}\Delta-8L^3_{h(\Delta)}\Delta^2\\&=\frac{4 \langle X_{k}^\Delta, g_\Delta(X_k^\Delta)\rangle^2}{|X_{k}^\Delta|^4}\Delta-o(\Delta),\endaligned$$
where
$$\aligned B:&=2\langle X_k^\Delta, f_\Delta(X_k^\Delta)\rangle\Delta+|f_\Delta(X_k^\Delta)|^2\Delta^2\\&\quad+|g_\Delta(X_k^\Delta)|^2(\Delta B_k)^2+2\langle f_\Delta(X_k^\Delta),g_\Delta(X_k^\Delta)\rangle\Delta\cdot \Delta B_k.\endaligned$$

Moreover, we have
$$\aligned\mathbb{E}(\xi^3_k|\mathscr{F}_{k\Delta})&\le C(L^3_{h(\Delta)}\Delta^3+L^6_{h(\Delta)}\Delta^6
+L^6_{h(\Delta)}\Delta^3\\&\quad+L^3_{h(\Delta)}\Delta^2+L^4_{h(\Delta)}\Delta^3
+L^4_{h(\Delta)}\Delta^2\\&\quad+L^5_{h(\Delta)}\Delta^4
+L^6_{h(\Delta)}\Delta^5+L^6_{h(\Delta)}\Delta^4)\\&=o(\Delta),\endaligned$$
where $C$ is a positive constant independent of $\Delta.$

Therefore,
$$\aligned\mathbb{E}(|X_{k+1}^\Delta|^p|\mathscr{F}_{k\Delta})&\le|X_{k}^\Delta|^p\left[1+\frac{p}{2}\Big(\frac{2\langle X_{k}^\Delta, f_\Delta(X_k^\Delta)\rangle+|g_\Delta(X_k^\Delta)|^2}{|X_{k}^\Delta|^2}\right.\\&
\qquad\left.+(p-2)\frac{\langle X_k^\Delta,g_\Delta(X_k^\Delta)\rangle^2}{|X_{k}^\Delta|^4}\Big)\Delta+o(\Delta)\right]\\&
\le|X_{k}^\Delta|^p(1
-p\lambda\Delta+o(\Delta)).\endaligned$$
We have use Lemma \ref{l2} in the last inequality.

Now for any given $\varepsilon\in(0,\lambda)$ sufficiently small, we can choose $\Delta^*\in(0,1)$ sufficiently small such that for all $\Delta\in(0,\Delta^*),$ it follows that $\frac{o(\Delta)}{\Delta}\le p\varepsilon$.

Thus,

$$\aligned\mathbb{E}(|X_{k+1}^\Delta|^p|\mathscr{F}_{k\Delta})&
\le|X_{k}^\Delta|^p(1
-p(\lambda-\varepsilon)\Delta).\endaligned$$

Taking expectation on both sides yields

$$\aligned\mathbb{E}(|X_{k+1}^\Delta|^p|)&
\le\mathbb{E}(|X_{k}^\Delta|^p)(1
-p(\lambda-\varepsilon)\Delta).\endaligned$$

Thus, for any $k\ge1,$
\begin{equation}\label{ju}\aligned\mathbb{E}(|X_{k}^\Delta|^p|)&
\le|x(0)|^p(1
-p(\lambda-\varepsilon)\Delta)^k\le |x(0)|^pe^{-kp(\lambda-\varepsilon)\Delta}.\endaligned\end{equation}

That is the MTEM method (\ref{num}) is $p$-th moment exponentially stable of order $p(\lambda-\varepsilon).$ We have proved (\ref{shuzhi}).

By using Chebyshev inequality and Borel-Cantelli lemma, (\ref{shuzhi1}) could be obtained in the same way as \cite{HMY}.

To prove (\ref{shuzhi2}), given any fixed $\Delta>0$ and suppose $k\Delta\le t<(k+1)\Delta.$ Then
$$x_\Delta(t)-\bar{x}_\Delta(t)=x_\Delta(t)-X_k^\Delta=f_\Delta(X_k^\Delta)(t-k\Delta)+g_\Delta(X_k^\Delta)(B(t)-B(k\Delta)).$$

Since $0<p<1$, we have
$$\mathbb{E}|x_\Delta(t)-\bar{x}_\Delta(t)|^p\le \Delta^p\mathbb{E}|f_\Delta(X_k^\Delta)|^p+\mathbb{E}(|g_\Delta(X_k^\Delta)|^p\mathbb{E}(|B(t)-B(k\Delta)|^p|\mathscr{F}_{k\Delta}).$$

Now $f_\Delta$ and $g_\Delta$ are globally Lipschitz continuous by Lemma \ref{l1}, and notice that $B(t)-B(k\Delta)$ is independent of $\mathscr{F}_{k\Delta},$ then
$$\aligned\mathbb{E}|x_\Delta(t)-\bar{x}_\Delta(t)|^p&\le \Delta^p\mathbb{E}(3L_{h(\Delta)}|X_k^\Delta|)^p
+\mathbb{E}(3L_{h(\Delta)}|X_k^\Delta|)^p\Delta^\frac{p}{2}\big)\\&
\le 3^pL_{h(\Delta)}^p(\Delta^p+\Delta^\frac{p}{2})\mathbb{E}(|X_k^\Delta|^p).\endaligned$$

Notice that for any $k\Delta\le t<(k+1)\Delta$, $\mathbb{E}|\bar{x}_\Delta(t)|^p=\mathbb{E}(|X_k^\Delta|^p).$ Therefore, by (\ref{tiaoj}), for $\Delta$ small enough, there exists $C>0$ such that
$$\aligned\mathbb{E}|x_\Delta(t)|^p&\le\mathbb{E}|x_\Delta(t)-\bar{x}_\Delta(t)|^p+\mathbb{E}|\bar{x}_\Delta(t)|^p\\&
\le(3^pL_{h(\Delta)}^p(\Delta^p+\Delta^\frac{p}{2})+1)\mathbb{E}(|X_k^\Delta|^p)\le C\mathbb{E}(|X_k^\Delta|^p).\endaligned$$

So by (\ref{ju}) we just obtained above, we have
$$\aligned\limsup_{t\to\infty}\frac{\log\mathbb{E}(|x_\Delta(t)|^p)}{t}&\le \limsup_{k\to\infty}\sup_{k\Delta\le t<(k+1)\Delta}\frac{\log\mathbb{E}(|x_\Delta(t)|^p)}{t}\\&
\le\limsup_{k\to\infty}\sup_{k\Delta\le t<(k+1)\Delta}\frac{\log C+\log\mathbb{E}(|X_k^\Delta|^p)}{t}\\&
\le \limsup_{k\to\infty}\frac{\log C+p\log|x_0|-kp(\lambda-\varepsilon)\Delta}{(k+1)\Delta}\\&
=p(\lambda-\varepsilon).\endaligned$$

This completes the proof. $\square$

\section{Example}

In this section, let us present an example to interpret our theory.

Let $d=1.$ Consider the following scalar SDE:
\begin{equation}\label{sde1}
dx(t)=(x(t)+x^3(t))dt+2\sqrt{x^4(t)+2x^2(t)}dB(t).
\end{equation}

In this equation, the drift term $f(x)=x+x^3$ and the diffusion term $g(x)=2\sqrt{x^4+2x^2}$.

Notice that neither $f$ nor $g$ is linear growing, meanwhile $f$ is not one-sided Lipschitz continuous. Then the exponential stability of classical Euler-type numerical approximations such as EM method, backward EM method or $\theta$-EM method could not be applied for this example.

It is obvious that both $f$ and $g$ are locally Lipschitz continuous (with $L_R=(1+3R^2)\vee(2+2R)$). Moreover, since in this case
\begin{equation}
\aligned\frac{\langle x,f(x)\rangle+\frac{1}{2}|g(x)|^2}{|x|^2}+\frac{p-2}{2}\frac{\langle x,g(x)\rangle^2}{|x|^4}&=\frac{\langle x,f(x)\rangle+\frac{p-1}{2}|g(x)|^2}{|x|^2}\\&
=\frac{(4p-3)x^2+(2p-1)x^4}{x^2},\endaligned
\end{equation}
then (\ref{tiaoj1}) holds for $p=\frac{1}{2}$ and $\lambda=1.$ Therefore, equation (\ref{sde1}) admits a unique strong global solution.

Let $h(\Delta)=\sqrt{\frac{\Delta^{-\frac{1}{5}}-1}{3}}$ for $\Delta<4^{-5}.$ Then $1<h(\Delta)\to\infty$ as $\Delta\to0,$ moreover, we have
$$L_{h(\Delta)}^4\Delta=(1+3h^2(\Delta))^4\Delta=\frac{1}{1+3h^2(\Delta)}=\Delta^{\frac{1}{5}}\to0\ \textrm{as}\ \Delta\to0.$$

Thus, condition (\ref{tiaoj}) also holds for given $h$. Then by Theorem \ref{stab}, we know that the $\frac{1}{2}$-th moment of the unique global solution of (\ref{sde}) is exponentially stable with order $\frac{1}{2}$. Moreover, for the numerical approximation MTEM (\ref{num}), we have that for any $\varepsilon<1,$
\begin{equation}\label{li1}\limsup_{k\to\infty}\frac{\log\mathbb{E}(|X_{k}^\Delta|^\frac{1}{2})}{k\Delta}\le-\frac{1}{2}(1-\varepsilon)\end{equation}
and
\begin{equation}\label{li2}\limsup_{k\to\infty}\frac{\log|X_{k}^\Delta|}{k\Delta}\le-(1-\varepsilon).\end{equation}

For the continuous-times MTEM approximation (\ref{num2}), we also have
\begin{equation}\label{li3}\limsup_{t\to\infty}\frac{\log\mathbb{E}(|x_\Delta(t)|^\frac{1}{2})}{t}\le-\frac{1}{2}(1-\varepsilon).\end{equation}

Thus, the MTEM replicates the exponential stability of the exact solution for the given SDE.

\end{document}